\documentclass[12pt]{amsart}
 \RequirePackage{amsmath}
 \RequirePackage{amssymb}
 \RequirePackage{amsthm}
 \RequirePackage{epsfig}
 \usepackage{palatino}

\def\a{\alpha}               \def\g{\gamma}
\def\d{\delta}       \def\De{{\Delta}}    
\def\la{\lambda}           
              
\def\t{\theta}             
\def\z{\zeta}                 

\def\C{{\mathbb C}}  \def\D{{\mathbb D}}
  \def\N{{\mathbb N}}
\def\R{{\mathbb R}}

\def\({\left(}       \def\){\right)}

\newtheorem{prop}{\sc Proposition}
\newtheorem{lem}{\sc Lemma}
\newtheorem{thm}{\sc Theorem}
\newtheorem{cor}{\sc Corollary}

\newtheorem{other}{\sc Theorem}              

\begin{document}
\title[Derivative of a Blaschke product and normal weights]{On Blaschke products with derivatives in Bergman spaces with normal weights}
\author[A. Aleman and D. Vukoti\'c]{Alexandru Aleman and Dragan Vukoti\'c}
\address{Department of Mathematics, Lund University, P.O. Box 118, Lund, Sweden}
\email{aleman@maths.lth.se}
\address{Departamento de Matem\'aticas and ICMAT CSIC-UAM-UC3M-UCM, M\'odulo C-XV, Universidad Aut\'onoma de Madrid, 28049 Madrid, Spain}
 \email{dragan.vukotic@uam.es}
 \thanks{The second author is supported by the Spanish
grant MTM2006-14449-C02-02 and also partially by the Thematic
Network MTM2006-26627-E, both from the Ministerio de Ciencia e Innovaci\'on, Spain, as well as by the European ESF Networking Programme (Network ``Harmonic and Complex Analysis and Its Applications'').}
\date{21 May, 2009}
\begin{abstract}
We generalize a well-known sufficient condition for interpolating sequences for the Hilbert Bergman spaces to other Bergman spaces with normal weights (as defined by Shields and Williams) and obtain new results regarding the membership of the derivative of a Blashke product or a general inner function in such spaces. We also apply duality techniques to obtain further results of this type and obtain new results about interpolating Blaschke products.
\end{abstract}
\maketitle
\section*{Introduction}
The conditions under which the derivative of a Blaschke product belongs to certain spaces of analytic functions in the disk is a subject with long history. Membership of the derivative in the Hardy spaces $H\sp{p}$ was investigated in \cite{Pr}, \cite{AC1}, \cite{A1}, and \cite{Ch} and a study of when it belongs to $B\sp{p}$, the Banach envelope of $H^p$ with $0<p<1$, can be found in \cite{AC1, AC2}. For an overview of such results until the early 1980s, we refer the reader to \cite{Cl}.
\par
Membership of the derivative of a Blaschke product in weighted Berg\-man spaces $L\sp{p,\a}_a$ was studied in \cite{A2} and \cite{Ki} and more recently in \cite{Ku}, \cite{GPV1, GPV2, GPV3}, \cite{GP1, GP2}, and \cite{FM}. Recall that for $\a>-1$ the \textit{standard weighted Bergman space\/} $L^{p,\a}_a$ (often also denoted by $A\sp{p}_\a$) is defined as the set of all analytic functions $f$ in the unit disk $\D$ for which
\[
 \|f\|_{L^{p,\a}_a}^p = (\a+1) \int_\D |f(z)|^p (1-|z|^2)^\a  dA(z)<\infty\,.
\]
When $\a=0$, we write simply $L^p_a$ and refer to it as the \textit{Bergman space}. A direct application of the Schwarz-Pick lemma shows that the derivative of any Blaschke product $B$ belongs to $\cap_{0<p<1} L^p_a$. Rudin \cite{R} proved that there are Blaschke products whose derivative does not belong to $L^1_a$ and Piranian \cite{Pi} gave an explicit example. Ahern \cite{A2} found a necessary and sufficient condition for the membership of $B\sp{\prime}$ in the standard weighted Bergman space $L^{p,\a}_a$ expressed in terms of $|B|$.
\par
Since Ahern's condition is non-trivial to check, it is desirable to have other conditions, either necessary or sufficient, for the membership of $B^\prime$ in $L_a^p$ or in more general weighted Bergman space $L_a^p(w)$, where by a \textit{weight\/} $w$ we mean any strictly positive measurable function $w$ in $\D$ which replaces the function $(\a+1) (1-|z|^2)^\a$ in the definition. Such criteria are often given in terms of the moduli of the zeros of $B$ or in terms of their location: belonging to some Stolz angle, some separation (discreteness) condition, etc. It is precisely this setting in which a study was carried out in \cite{Ki}, \cite{GP1, GP2}, \cite{GPV1, GPV2, GPV3}, or \cite{FM}.
\par
Recall that a sequence $\{a_n\}$ in $\D $ is said to be \emph{separated\/} (or {\emph{uniformly discrete\/}) with \textit{constant of separation\/} $\d$ if
\begin{equation}\label{separated}
 \inf_{n\in\N} \prod_{k\neq n} \left| \frac{a_k-a_n}{1-\overline{a}_k
 a_n} \right| =\d > 0\,.
\end{equation}
By a well-known theorem of Carleson, this is equivalent to $(a_n)_{n=1}^\infty$ being an interpolating sequence for the space $H^\infty$ of bounded analytic functions. An important class of separated sequences is that of exponential sequences; it is known that, unlike in Rudin's example, the derivative of every Blaschke product whose zeros form an exponential sequence belongs to $L^1_a$.
\par
The main purpose of this article is to extend several earlier results to the context of weighted Bergman spaces with normal weights as defined by Shields and Williams \cite{SW}. For a non-negative radial weight on the unit disk $\D$, its values in $\D$ are completely determined by the values in the interval $[0,1)$. A positive function $w(r)$ in this interval is normal if, roughly speaking, it grows at the rate controlled by two positive powers of $1-r$ (see \cite{SW} and \cite{AS}, for example). As a ``rule of thumb'', a space with such weights behaves to a certain extent like the unweighted one.
\par
In Section~\ref{sec-normal} of the paper we review the necessary background and prove two results. One is a basic proposition of Shields-Williams type which explains the asymptotic behavior of certain weighted integrals and which is needed in the rest. The other is a sufficient condition for interpolation in weighted Bergman spaces. Our main tool in this section is the use of the so-called B\'ekoll\'e condition \cite{Be}.
\par
In Section~\ref{sec-main} we consider some conditions which are either necessary or sufficient for membership of the derivative of a Blaschke product in Bergman spaces with normal weights, thus generalizing several results obtained previously. We also study the integrability of the derivative of a general inner function in the disk in terms of the properties of its generalized counting functions. Our results allow us to deduce as a corollary one of the main results from \cite{FM}. Also,   several propositions from \cite{GPV1} follow for free as well. It follows from our  main theorems that there are interpolating Blaschke products which are in no $L_a^p$, $p>1$. An explicit example was given in \cite[Theorem~3.8]{GP1}. We remark that Pel\'aez \cite{Pe}, inspired by Piranian's construction \cite{Pi}, has recently constructed such a product which is not in $L_a^1$. It is not clear whether our method  could allow us to produce such examples, so there are some limitations to the techniques used here.
\par
In Section~\ref{sec-duality} we use Luecking's duality theorem for weighted Bergman spaces to obtain some further estimates on the integral of $|B^\prime|$. As a main result, our Theorem~\ref{cohn-an} provides the weighted Bergman space analogue of a theorem due to Cohn \cite{Ch} which relates the membership of $|B^\prime|$ in certain Hardy spaces with the orthogonal complement of the invariant subspace of $H^2$ generated by $B$.
\par
Section~\ref{sec-higher} includes some closing remarks. In it, we  mention some possible formulations of our results in terms of derivatives of higher order.

\section{Normal weights}
 \label{sec-normal}
We will typically consider measurable functions $w\,\colon\,\D\to [0,+\infty)$. Such a function is said to be \textit{radial\/} if $w(z)=w(|z|)$ for all $z$ in $\D$; obviously, its values are then  completely determined by the values on the interval $[0,1)$. Since such functions are usually important examples of weights for Bergman spaces, further special growth or integrability conditions are often imposed on them, especially in relation to certain integral operators; see \cite{AS}, for example. In this paper we will also consider one special and important class of radial weights: the so-called normal ones. Following Shields and Williams \cite{SW}, we will say that the function $w\,\colon\, [0,1)\to [0,\infty)$ is \textit{normal\/} if there exist real numbers $a$ and $b$ and $r_0\in (0,1)$ for which the two following conditions are fulfilled:
\[
 \frac{w(r)}{(1-r)^a} \nearrow \infty \quad \textrm{for}\quad
 r>r_0
\]
and
\[
 \frac{w(r)}{(1-r)^b} \searrow 0 \quad \textrm{for}\quad
 r>r_0 \,.
\]
One of the basic properties of radial normal functions is the following fact.
\par
\begin{lem}
 \label{lem-norm-w}
For each $t\in (0,1)$ there exist constants $c$ and $C$ depending
only on $t$ such that
\[
 c\,w(\z)\le w(z)\le C\,w(\z)
\]
whenever $|z|>r_0$ and $|z-\z|<t\,(1-|z|)$.
\end{lem}
\par
\begin{proof}
The proof is simple, so to illustrate the method it suffices to prove  only one of the inequalities, say in the case $|\z|\ge |z|$. Under this assumption, the obvious inequalities
\[
 \frac{w(|z|+t(1-|z|))}{(1-t)^b (1-|z|)^b} = \frac{w(|z|+t(1-|z|))}{(1-|z| - t(1-|z|))^b} \le \frac{w(|z|)}{(1-|z|)^b} = \frac{w(z)}{(1-|z|)^b}
\]
yield
\[
 w(\z)=w(|\z|)\le w(|z|+t(1-|z|))\le (1-t)^b w(z)
\]
whenever $|z-\z|<t\,(1-|z|)$.
\end{proof}
\par
Given a normal function $w$, denote by $a_w$ and $b_w$ the optimal
indices
\[
 a_w=\inf\{a\,\colon\,\frac{w(r)}{(1-r)^a}\nearrow\infty\,, r>r_a\}
 \,, \quad b_w=\sup\{b\,\colon\,\frac{w(r)}{(1-r)^b} \searrow 0\,,
 r>r_b\}\,.
\]
Note that $a_w\geq b_w$ and also that it is not required that the
functions
\[
 \frac{w(r)}{(1-r)^{a_w}} \,, \qquad w(r) (1-r)^{b_w}
\]
be increasing and decreasing respectively.
\par
We list below some estimates for normal functions. They should be
compared with the usual integral estimates for the standard radial
weights. For example, it is well known that
\[
 I_{\g,m}(\la)\asymp\int_\D \frac{(1-|z|)^\g}{|1-\overline{\la} z|^{m+1}}
 \,dA(z) = O \( (1-|\la|)^{\g-m+1} \)
\]
as $|\la|\to 1$ (see \cite[Theorem~1.7]{HKZ}, for example). As is usual, throughout the paper the notation $u(\la)\asymp v(\la)$ will mean that the quotient of two positive functions $u$ and $v$ is bounded from above and below as $|\la|\to 1$.
\par
The estimates below are essentially known to the experts and can be
attributed to Shields and Williams. Certainly, some lemmas in their paper \cite{SW} have a similar flavor and use similar proofs.
\par
\begin{prop} \label{thm-sh-wil}
Let $w$ be a normal function  and let $m\in \R$. For $\la\in \D$,
let
\begin{align*}
 &I_{m,w}(\la)=\int_\D \frac{w(|z|)}{|1-\overline{\la}
 z|^{m+2}} d A(z) \,,\\&
 J_{m,w}(\la)=\int_\D\log\left|\frac{1-\overline{\la}z}{\la-z}\right|
 (1-|z|)^{-m-2}w(|z|)dA(z)\,.
\end{align*}
\begin{enumerate}
 \item[(i)] If $m>-1$, $a_w<m$, and $b_w>-1$ then
\[
 I_{m,w}(\la)\asymp w(|\la|) (1-|\la|)^{-m}  \,.
\]
\item[(ii)] If  $m\in \R$, $a_w<m+1$, and $b_w>m$ then
\[
 J_{m,w}(\la)\asymp w(|\la|) (1-|\la|)^{-m} \,.
\]
\end{enumerate}
Assume also that:
\par
(*)  There exists $\a$ such that $0<\a<1$ and the function
\[
 \frac{w(r)}{(1-r)^{a_w}} \log^\a \frac1{1-r}
\]
is increasing for $r>r_0$.
\par
Then:
\begin{enumerate}
\item[(iii)] If $m>-1$, $a_w=m$ and $b_w>-1$
\[
 I_{m,w}(\la)=O\( w(|\la|) (1-|\la|)^{-m} \log \frac1{1-|\la|} \)
 \quad when \quad |\la|\to 1\,.
\]
\item[(iv)] If $m\in \R$,  $a_w=m+1$, and $b_w>m$ then
\[
 J_{m,w}(\la)=O\( w(|\la|) (1-|\la|)^{-m}\log\frac1{1-|\la|} \) \quad as \quad
 |\la|\to 1\,.
\]
\end{enumerate}
\end{prop}
\par
\begin{proof}
We only prove parts (i) and (ii) in order to illustrate the basic
technique.
 Let $a$, $b$ be arbitrary real numbers such that $m>a>a_w$ and
$-1<b<b_w$. \\ (i) By applying the usual splitting argument for the
normal weights and taking into account that the function $(1-r)^{-a} w(r)$ is increasing on $(r_0,1)$ and $(1-r)^{-b} w(r)$ is decreasing on $(r_0,1)$, we obtain for $r_0<|\la|<1$
\begin{eqnarray*}
 I_{m,w}(\la) &=& \int_\D \frac{w(z)}{|1-\overline{\la}
 z|^{m+2}} d A(z)
 =\int_0^1 r w(r)\( \int_0^{2\pi} \frac{d\t}{|1-
 \overline{\la} z|^{m+2}} \) \frac{d r}{\pi}
\\
 &\le & C \int_{r_0}^1 \frac{r w(r)}{(1-|\la| r)^{m+1}} d r
\\
 &\le & C \( \int_{r_0}^{|\la|} \frac{w(r)}{(1-|\la| r)^{m+1}}
 d r + \int_{|\la|}^1 \frac{w(r)}{(1-|\la| r)^{m+1}} d r \)
\\
 &\le & C  (1-|\la|)^{-a} w(|\la|) \int_{r_0}^{|\la|}
 \frac{(1-r)^{a}}{(1-|\la| r)^{m+1}} d r
\\
 & & +  C (1-|\la|)^{-b} w(|\la|) \int_{|\la|}^1
 \frac{(1-r)^{b}}{(1-|\la| r)^{m+1}} d r
\\
 &\le & C  (1-|\la|)^{-a} w(|\la|) \int_{0}^{|\la|}
 (1-r)^{a-m-1}
\\
 & & +(1-|\la|)^{-b-m-1} w(|\la|) \int_{|\la|}^1(1-r)^{b} d r
\\
 &\le & C w(|\la|) (1-|\la|)^{-m} \,.
\end{eqnarray*}
The constant $C$ above may differ from one appearance to another.
Similarly, since $-1<b<b_w\leq a_w<m$,
\begin{align*}
 I_{m,w}(\la)&\geq C\int_{r_0}^{|\la|}\frac{w(r)}{(1-|\la| r)^{m+1}}
 d r\geq C(1-|\la|)^{-b}w(|\la|)\int_{r_0}^{|\la|}
 \frac{(1-r)^b}{(1-|\la| r)^{m+1}}\\&
 \geq
 C(1-|\la|)^{-m}w(|\la|)\int_1^{\frac{1-r_0}{1-|\la|}}
 \frac{x^b}{(1+x)^{m+1}}dx\geq C(1-|\la|)^{-m}w(|\la|)\,.
\end{align*}
(ii)  By Jensen's formula we have
$$
J_{m,w}(\la)=\log\frac1{|\la|}\int_0^{|\la|} w(r)(1-r)^{-m-2}rdr +
\int_{|\la|}^1\log\frac1{r} w(r)(1-r)^{-m-2}rdr\,,
$$
and the estimates for $J_{m,w}(\la)$ follow with the same argument as above.
\end{proof}
If the normal weight $w$ is integrable, for example, if $b_w>-1$,  we can consider the weighted Bergman space $L_a^p(w),~p>0$, which consists of analytic functions $f$ in $\D$ with
$$
 \|f\|_{p,w}^p=\int_\D|f|^pwdA<\infty\,.
$$
Using Lemma \ref{lem-norm-w}, it follows easily that if $f\in
L_a^p(w)$ then
$$
 |f(z)|\leq C\(w(|z|)(1-|z|)^2\)^{-1/p}\|f\|_{p,w}\,,
$$
for some constant $C>0$ and all $z\in \D$. By Proposition
\ref{thm-sh-wil} (i) we have that
$$
 \sup_{\|f\|_{p,w}\le 1}|f(z)|\asymp \(w(|z|)(1-|z|)^2\)^{-1/p}\,,
$$
if $b_w>-1$. We can then consider interpolating sequences for these
spaces. A sequence $(z_n)_{n=1}^\infty$ in $\D$ is called an \textit{interpolating
sequence\/} for $L_a^p(w)$ if for every $(a_n)\in \ell^p$ there exists
$f\in L_a^p(w)$ such that
$$
 f(z_n) \(w(|z_n|)(1-|z_n|)^2\)^{1/p}=a_n\,,
$$
for all $n$. If $p=2$, it is well known that interpolating sequences for Hardy spaces are interpolating for $L_{a,w}^2$ as well. (For example, this follows from \cite[p. 23]{Seip} by recalling the fact which is a ```folk knowledge'': the pointwise multipliers of $L_{a,w}^2$ into itself are precisely the bounded analytic functions.) This interpolation result can be extended for the other values of $p>1$.

\begin{thm}\label{interpolation} If $b_w>-1$, then every interpolating sequence for $H^\infty$ is
interpolating for $L_a^p(w),~p> 1$.
\end{thm}

\begin{proof} Our proof uses B\'ekoll\'e's theorem \cite{Be} and
Luecking's duality theorem \cite{L}. To be more explicit, if
$\frac{1}{p}+\frac1{q}=1$, let $\g>\frac{a_w}{p}-\frac{1}{q}$ and
use a similar argument to the above to obtain that
$$\(\int_r^1w(t)dt\)\(\int_r^1(1-t)^{\g q}w^{-q/p}(t)dt\)^{p/q}\leq
C(1-r)^{p\g+p}\,,$$ which implies that $w(|z|)(1-|z|)^{-\gamma}$
belongs to the B\'ekoll\'e class $B_p(\gamma)$, or equivalently,
$v(z)=w^{-q/p}(z)(1-|z|)^{(q-1)\g}$ belongs to $B_q(\gamma)$. Then
by a result in \cite{Be} it follows that the sublinear
operator $P_\g$ defined by \begin{equation}\label{projection}P_\g
f(z) =\int_\D\frac{(1-|\z|)^{\g}}{|1-\overline{\z}z|^{\g+2}}
|f(\z)|dA(\z)\,,\end{equation} is bounded from
$L^q((1-|z|)^{q\g}w^{-q/p}dA)$ into itself. The duality theorem in
\cite{L} then identifies the dual of $L_a^p(w)$ with $L_a^q(v)$ via
the standard weighted Bergman pairing with index $\g$. In other
words, the norm in $L_a^p(w)$ is equivalent to
\begin{equation}\label{duality}
 \sup \left\{ \left| \int_\D f(z) \overline{g(z)}(1-|z|^2)^\g d A
 \right| \,\colon\, g\in
 L_a^q(w)\,, \|g\|_{q,v}\le 1\right\} \,.
\end{equation}
Now let $(z_n)_{n=1}^\infty$ be an interpolating sequence for $H^\infty$,  let
$B$ be the Blaschke product with zeros $z_1,\ldots, z_n,\ldots$ and
let $a=(a_n)\in \ell^p$. We claim that the solution $f_a$ of the
interpolation problem defined above is given by
$$f_a(z)=\sum_na_n(w(|z_n|)(1-|z_n|)^2)^{-1/p}
\frac{(1-|z_n|^2)^{\g+1}B(z)}{(1-\overline{z}_nz)^{\g+1}(z-z_n)B'(z_n)}
\,.$$ In order to verify this, it suffices to show that
$$
 \|f_a\|_{p,w}\leq C\|a\|_{\ell^p}\,,
$$
for some constant $C$ and all sequences $a\in \ell^p$ with finitely many nonzero terms. Denote by $B_n$ the Blaschke product with zeros $\{z_k:~k\neq n\}$ taking multiplicities into account, and recall that the values $|B_n(z_n)|$ are bounded away from zero by \eqref{separated}. By
 \eqref{duality} we have
\begin{align*}
 &\|f_a\|_{p,w}\leq\\& \sup_{\|g\|_{q,v}\le 1}\sum_n
 \frac{|a_n|}{|B_n(z_n)|}w^{-1/p}(|z_n|)
 (1-|z_n|)^{\g+2/q}\int_\D\frac{|g(z)||B_n(z)|(1-|z|)^\g}{|1-
 \overline{z}_nz|^{\g+2}}dA(z)\\& \le C\sup_{\|g\|_{q,v}\le
 1}\sum_n|a_n|w^{-1/p}(|z_n|) (1-|z_n|)^{\g+2/q}P_\g g(z_n)\\&\le C
 \|a\|_{\ell^p}\sup_{\|g\|_{q,v}\le 1}\(\sum_nw^{-q/p}(|z_n|)
 (1-|z_n|)^{q\g+2}P^q_\g g(z_n)\)^{1/q}\,.
\end{align*}
Since $(z_n)_{n=1}^\infty$ is separated, there is a fixed $R\in (0,1)$ such that the disks
\[
 \De_n = \{\zeta\in\D\,\colon\,|\zeta-z_n|<R\,(1-|z_n|)\}\,,
\]
are pairwise disjoint. Moreover, it is easy to show that $P_\g g$ is
almost constant on these disks, \textit{i.e.}, there exists $c_0>0$
independent of $g$ such that
$$
 c_0^{-1}P_\g g(z)\le P_\g g(\z)\le c_0P_\g g(z)\,, \quad z,\z\in
\De_n\,.
$$
This, together with Lemma~\ref{lem-norm-w}, leads to
\begin{align*}
 \sup_{\|g\|_{q,v}\le 1}\sum_n&w^{-q/p}(|z_n|)
 (1-|z_n|)^{q\g+2}P^q_\g g(z_n)\\&\leq C\sup_{\|g\|_{q,v}\le
 1}\sum_nw^{-q/p}(|z_n|)(1-|z_n|)^{q\g}\int_{\De_n}P^q_\g
 g(z)dA(z)\\& \le C\sup_{\|g\|_{q,v}\le 1}\sum_n\int_{\De_n}P^q_\g
 g(z)w^{-q/p}(z)(1-|z|)^{q\g}dA(z)\\&\le \sup_{\|g\|_{q,v}\le 1}
 \int_{\cup_n\De_n}P^q_\g g(z)w^{-q/p}(z)(1-|z|)^{q\g}dA(z)\,.
\end{align*}
and the result follows by  \eqref{projection}.
\end{proof}
\section{Main results}
 \label{sec-main}
\par
In what follows, for a nonconstant analytic function
$f$ in $\D$ and a function $u:\D\to \R_+$ we will consider the sum
$$
 \sum_{f(z)=\zeta}u(z)\,,
$$
the multiplicities being taken into account. In particular, such sums will be useful to denote the summation over all the zeros of a Blaschke product. We remind the reader that most results obtained in the literature on the integrability of $B^\prime$ are obtained under certain special assumptions on the distribution or moduli of the zeros of $B$.
\par
Our first statement in this section generalizes or adds new information to several known results. Besides complementing the well-known statement (attributed to Rudin) that
\[
 \|B^\prime\|_{1}\le c \sum_{B(z)=0} (1-|z|)\log\frac1{1-|z|}\,,
\]
it also includes as special cases several propositions from the recent papers  such as \cite{GPV1} and \cite{FM}.
\begin{thm}\label{main1}
Let $w$ be a normal function and let $1/2<p<\infty$.
\begin{enumerate}
 \item[(i)] Assume that   $a_w<2 p-2,~b_w>-1$, if  $1/2<p\le 1$,
 and that $a_w<p-1,~b_w>p-2$, if  $p>1$.  Then there
 exists a positive constant $c_{p,w}$ such that for every Blaschke
 product $B$ we have
\begin{equation}
 \|B^\prime\|^p_{p,w}\le c_{p,w} \sum_{B(z)=0} (1-|z|)^{2-p}
 \,w(|z|)\,.
 \label{ibp-upper}
\end{equation}
  If \/$1/2<p\le 1$ and  $a_w=2 p-2,~b_w>-1$, or if $p>1$,
 $a_w=p-1,~b_w>p-2$, and   $w$ satisfies
the condition (*) in Proposition \ref{thm-sh-wil}, then there
 exists a positive constant $c_{p,w}$ such that for every Blaschke
 product $B$ we have
\begin{equation}
 \|B^\prime\|^p_{p,w}\le c_{p,w} \sum_{B(z)=0}
 (1-|z|)^{2-p}\,w(|z|)\log\frac1{1-|z|}
 \,.
 \label{ibp-upper-lim}
\end{equation}
\item[(ii)] Suppose that the zero set  of the
Blaschke product $B$ is separated with separation constant
$\delta>0$. Suppose that either $a_w<2p-2$, $b_w>-1$, when $1/2<p\le 1$, or $a_w<p-2$, $b_w>-1$, when $p>1$. Then, in both cases, there exists a positive constant $c_{p,w,\d}$ such that
 \begin{equation}
  \sum_{B(z)=0}^\infty (1-|z|)^{2-p}
 \,w(|z|)\le c_{p,w,\delta}\|B^\prime\|^p_{p,w}\,.
 \label{ibp-lower}
\end{equation}
If either $1/2<p\le 1$ and $a_w=2p-2$, $b_w>-1$, or $p>1$,
$a_w=p-1$, $b_w>p-2$, and $w$ also satisfies the condition (*) in Proposition \ref{thm-sh-wil}, then there exists a positive constant $c_{p,w,\delta}$ such that for every Blaschke product $B$ we have
 \begin{equation}
  \sum_{B(z)=0}^\infty (1-|z|)^{2-p}
 \,w(|z|)\le c_{p,w,\delta}\int_\D|B'(z)|^pw(|z|)
 \log\frac1{1-|z|}dA(z)\,.
 \label{ibp-lower-limit}
\end{equation}
 \end{enumerate}
\end{thm}
\begin{proof} To prove (i) we assume first that $1/2<p\le 1$. Start off with the inequality
\begin{equation*}
 |B^\prime(\zeta)|\le \sum_{B(z)=0} \frac{1-|z|^2}{|1-\overline{z}
 \z|^2} \,,
 \label{deriv-upper}
\end{equation*}
which is quite straightforward to deduce (see \cite{FM}, for example). For $p\le 1$ this implies
\[
 \int_\D |B^\prime|^p\,w\,d A\le \sum_{B(z)=0} (1-|z|^2)^p \int_\D
 \frac{w(\zeta)}{|1-\overline{z} \zeta|^{2 p}}\,d A(z) \,.
\]
Thus for $1/2<p\le 1$, part  (i) follows by a direct application of
Theorem~\ref{thm-sh-wil} (i) and (iii), with $m=2p-2$.
\par
In the case $p>1$ we can apply the Schwarz-Pick lemma together with the elementary inequality $(1-x^2)^p\leq -2\log x$, $x\in [0,1]$, to obtain
\begin{align*}
 \int_\D |B^\prime|^p\,w\,d A &\le
 \int_\D \frac{(1-|B(\z)|^2)^p}{(1-|\z|^2)^p}\,w(|\z|)\,d A(\z)\\&
 \le -2\int_\D \log |B(\z)|\,w(|\z|)(1-|\z|^2)^{-p}\,d A(\z)\\&=
 2\sum_{B(z)=0}\int_\D\log\left|\frac{1-\overline{z}\z}{z-\z}\right|
 (1-|\z|)^{-p}w(|\z|)dA(\z)\\&
 =2\sum_{B(z)=0}J_{p-2,w}(z)\,,
 \end{align*}
and the result follows by Theorem~\ref{thm-sh-wil} (ii) and (iv).
\par
(ii) Assume that the zero set of $B$ is separated with the constant of
separation $\delta>0$. Then there is a fixed $R\in (0,1)$ depending
only on $\delta$  such that the disks
\[
 \De_z = \{\zeta\in\D\,\colon\,|\zeta-z|<R\,(1-|z|)\}\,,\quad
 B(z)=0\,,
\]
are pairwise disjoint. Moreover,
$$
 \rho(R)=\sup\{|B(\z)|\,\colon\, \z\in \cup_{B(z)=0}\De_z\}<1\,.
$$
See, for example, \cite[pp. 681--682]{GPV1} or \cite{N}.
\par
By Lemma \ref{lem-norm-w}, the obvious estimate for the area of the disk $\De_z$, and an application of the standard Green's formula (see, \textit{e.g.}, \cite[p.~16]{DS} for a formulation in the  smooth case),  and a routine Laplacian computations, for every $\eta\geq 0$ we have
\begin{align*}
 (1-&\rho(R)^2)^\eta\sum_{B(z)=0} (1-|z|)^{2-p}w(|z|)\\
&\leq C \sum_{B(z)=0} \int_{\De_z}(1-|B(\z)|^2)^\eta
 w(|\z|)(1-|\z|)^{-p}dA(\z)\\
&\leq C\int_\D(1-|B(\z)|^2)^\eta w(|\z|)(1-|\z|)^{-p}dA(\z)\\
& =-\frac{C}{2\pi} \int_\D\Delta  (1-|B(\xi)|^2)^\eta \int_\D\log\left|\frac{1-\overline{\z}\xi}{\z-\xi}\right|
 w(|\z|)(1-|\z|)^{-p}dA(\z)dA(\xi)\\
&\leq \frac{2\eta^2 C}{\pi} \int_\D
(1-|B(\xi)|^2)^{\eta-2}|B'(\xi)|^2 J_{p-2,w}(\xi)dA(\xi) \\& \leq
 \frac{2\eta^2 C K}{\pi} \int_\D (1-|B(\xi)|^2)^{\eta-2}|B^\prime(\xi)|^2 (1-|\xi|)^{2-p}w(|\xi|)dA(\xi)\,,
\end{align*}
where the last inequality follows from Proposition~\ref{thm-sh-wil}~(ii). Now if $p\leq 2$, we choose $\eta=2$ and apply again the
Schwarz-Pick lemma to obtain
$$
 (1-\rho(R)^2)^\eta\sum_{B(z)=0} (1-|z|)^{2-p}w(|z|)\leq
 \frac{8 C K}{\pi} \int_\D |B^\prime(\xi)|^p w(|\xi|)dA(\xi)\,.
 $$
In the case when $p>2$, we choose $\eta=2+\frac{p-2}{p}$ and apply H\"older's inequality (in addition to an argument used earlier in the proof of part (i) and Proposition~\ref{thm-sh-wil}~(ii)) to obtain from above
\begin{align*}
 (1-\rho(R)^2)^\eta&\sum_{B(z)=0} (1-|z|)^{2-p}w(|z|) \\
&\leq C_p \|B^\prime\|^{\frac2{p}}_{p,w}
 \left(\int_\D(1-|B(\xi)|^2)(1-|\xi|)^{-p}w(|\xi|)dA(\xi)
 \right)^{\frac{p-2}{p}} \\
& \leq C'_p\|B^\prime\|^{\frac2{p}}_{p,w}
 \left(-2 \int_\D\log|B(\xi)|(1-|\xi|)^{-p}w(|\xi|)dA(\xi)
 \right)^{\frac{p-2}{p}}\\
& \le C_p^{\prime\prime}\|B^\prime\|^{\frac2{p}}_{p,w}
 \( \sum_{B(z)=0}J_{p-2,w}(z)\)^{\frac{p-2}{p}} \\
& \le C_p^{\prime\prime\prime}\|B^\prime\|^{\frac2{p}}_{p,w}\left(\sum_{B(z)=0}
 (1-|z|)^{2-p}w(|z|)\right)^{\frac{p-2}{p}}\,,
\end{align*}
and the first estimate follows. The proof of the second inequality
in (ii) is almost identical and will be omitted.
\end{proof}
Let us observe right away that one of the main results of Fricain and Mashregi \cite{FM} is a direct consequence of Theorem \ref{main1}~(i). In
fact, our result together with the obvious estimates
\begin{align*}
 \frac1{\varepsilon}\int_{|z|<1/2}
 |B^\prime(z)|^pw(|z|)(1-|z|)^\varepsilon dA(z)& \leq
 \sup_{1/2<r<1}w(r)(1-r)\int_0^{2\pi} |B^\prime(re^{it})|^p dt\\
& \leq
 2\int_{|z|<1/2}|B^\prime(z)|^pw(|z|)dA(z)\,,
\end{align*}
where $\varepsilon>0$, lead to the following result.
\begin{cor}\label{mashregi}
If $w$ and $p$ satisfy the conditions in Theorem \ref{main1}~(i)
then there exists a constant $c_{p,w}>0$ such that
$$
 \sup_{1/2<r<1}w(r)(1-r)\int_0^{2\pi} |B^\prime(re^{it})|^p dt
 \leq \sum_{B(z)=0} (1-|z|)^{2-p} \,w(|z|)\,,
 $$
 when $a_w<2p-2,~1/2<p\le 1$, or $a_w<p-1,~p>1$, and
$$
 \sup_{1/2<r<1}w(r)(1-r)\int_0^{2\pi} |B^\prime(re^{it})|^p dt
 \leq \sum_{B(z)=0} (1-|z|)^{2-p}\,w(|z|)\log\frac{1}{1-|z|}\,,
$$
 when $a_w=2p-2,~1/2<p\le 1$, or $a_w=p-1,~p>1$.
\par
If $w$ and $p$ satisfy the conditions in Theorem \ref{main1} (ii) and the zero set of the Blaschke product $B$ is separated with separation constant $\d>0$, then for every $\varepsilon>0$ there exists a positive constant $c_{p,w,\delta,\varepsilon}$ such that
$$
  \sum_{n=1}^\infty (1-|z|)^{2-p+\varepsilon}
 \,w(|z|)\le c_{p,w,\delta,\varepsilon} \sup_{1/2<r<1}w(r)(1-r)
 \int_0^{2\pi} |B^\prime(re^{it})|dt\,.
$$
\end{cor}
\par
The lower estimate given in the second part of Theorem \ref{main1}
cannot be expected to hold for general Blaschke products. For
example, if $p=1$ and $B$ is a finite Blaschke product then it can be seen that $\|(B^n)^\prime\|_{1,w}=o(n)$, as $n\to \infty$, while
$$
 \sum_{B^n(z)=0}(1-|z|)w(|z|)=n\sum_{B(z)=0}(1-|z|)w(|z|)\,,
$$
for any normal weight $w$. From this point of view, the limit cases
$a_w=2p-2,~1/2<p\le 1$, and $a_w=p-1,~p>1$, raise even more delicate
problems.
\par
It turns out that an asymptotically sharp estimate  of the integrals
considered here can be obtained even for general inner functions,
if we consider pre-images of several points in the disk. To this end
we introduce, for a normal weight $w$, $p>0$ and analytic functions
$f$ in $\D$, the \textit{generalized counting functions}
\begin{equation}\label{counting1}
N_{f,w,p}(\z)=\sum_{f(z)=\z}(1-|z|)^{2-p}w(|z|)\,,
\end{equation}
and
\begin{equation}\label{counting2}
N^\ell_{f,w,p}(\z)=\sum_{f(z)=\z}(1-|z|)^{2-p}w(|z|)\log\frac1{1-|z|}\,.
\end{equation}
This notation allows us to formulate a result on integrability of the derivative of a general inner function rather than just for a Blaschke product.
\begin{thm}\label{inner-est}
Let $w$ be a normal function and let $1/2<p<\infty$.  Then there
 exists a positive constant $c_{p,w}$ such that for every
inner function $\theta$ and every $\varepsilon\in (0,1/2)$ we have
\begin{equation}
c^{-1}_{p,w} \int_{|\z|<\varepsilon}N_{\theta,p,w}(\z)dA(\z)\le
\|\theta^\prime\|^p_{p,w}\le \varepsilon^{-2}c_{p,w}
\int_{|\z|<\varepsilon}N_{\theta,p,w}(\z)dA(\z) \,,
 \label{counting-est}
\end{equation}
if $1/2<p\le 1$, and $a_w<2 p-2,~b_w>-1$, or $p>1$ and
$a_w<p-1,~b_w>p-2$. Moreover, in the limit case, when $1/2<p\le1,~
a_w=2 p-2,~b_w>-1$, or  $p>1, a_w=p-1,~b_w>p-2$, if
    $w$ satisfies
the condition (*) in Proposition \ref{thm-sh-wil} then
\begin{equation}
\|\theta^\prime\|^p_{p,w}\le \varepsilon^{-2}c_{p,w}
\int_{|\z|<\varepsilon}N^\ell_{\theta,p,w}(\z)dA(\z) \,,
 \label{counting-est1}
\end{equation}
and
\begin{equation}
 \int_{|\z|<\varepsilon}N_{\theta,p,w}(\z)dA(\z) \le c_{p,w}
\int_\D|\theta'(\z)|^pw(|\z|)\log\frac{1}{1-|\z|}dA(\z)\,.
 \label{counting-est2}
\end{equation}
\end{thm}
\begin{proof} By Frostman's theorem \cite[Theorem~6.4]{G}, for almost every $\xi\in\D$, the function $\theta_\xi= (\theta-\xi) (1-\overline{\xi}\theta)^{-1}$ is a Blaschke product. If $|\xi|<\varepsilon<1/2$ we can apply Theorem
\ref{main1} (i) to obtain
$$
 \|\theta'\|_{p,w}^p\leq \frac{C}{\varepsilon^2}\int_{|\xi|<\varepsilon}
 \|\theta'_\xi\|^p_{p,w}dA(\xi)\leq \frac{c_{p,w}}{\varepsilon^2}
 \int_{|\xi|<\varepsilon}\(\sum_{\theta_\xi(z)=0}(1-|z|)^{2-p}w(|z|)\)
 dA(\xi)\,,
$$
and the upper estimate in (\ref{counting-est}) follows. The proof of (\ref{counting-est1}) follows exactly the same steps, using Theorem \ref{main1}~(iii).
\par
To get the lower estimate in (\ref{counting-est}), assume first that $p\leq 2$, and use the standard change of variable formula \cite[p.~186]{Sh} together with the Schwarz-Pick lemma to obtain
\begin{align*}
\int_{|\z|<\varepsilon} N_{\theta,p,w}(\z)dA(\z)
&=\int_{|\theta|<\varepsilon}|\theta'(\z)|^2w(|\z|) (1-|\z|)^{2-p}
dA(\z)\\&\le \int_\D|\theta'(\z)|^pw(|\z|)dA(\z)\,.
\end{align*}
If
$p>2$ we have as above
\begin{align*}
&C\varepsilon^2\int_{|\z|<\varepsilon}N_{\theta,p,w}(\z)dA(\z)\\&\le
\int_{|\xi|<\varepsilon}\int_{|\z|<\varepsilon}
\(\frac{(1-|\xi|^2)(1-|\z|^2)}{|1-\overline{\xi}\z|^2}\)^{p-2}
N_{\theta,p,w}(\z)dA(\z)dA(\xi)\\&
=\int_{|\xi|<\varepsilon}\int_{|\theta|<\varepsilon}(1-|\theta_\xi(\z)|^2)^{p-2}
|\theta'(\z)|^2w(|\z|) (1-|\z|)^{2-p} dA(\z)dA(\xi) \,.
\end{align*}
By H\"older's inequality and Proposition \ref{thm-sh-wil} (ii) we
obtain
\begin{align*}&C\varepsilon^2\int_{|\z|<\varepsilon}N_{\theta,p,w}(\z)dA(\z)
\\&\le \pi \varepsilon^{p}\|\theta'\|_{p,w}^{2}
\(\int_{|\xi|<\varepsilon}\int_\D\frac{(1-|\theta_\xi(\z)|^2)^{p}}{(1-|\z|)^p}
w(|\z|)dA(\z)dA(\xi)\)^{\frac{p-2}{p}}\\& \leq \pi
\varepsilon^p\|\theta'\|_{p,w}^2
\(-2\int_{|\xi|<\varepsilon}\int_\D\log|\theta_\xi|(1-|\z|)^{-p}
w(|\z|)dA(\z)dA(\xi)\)^{\frac{p-2}{p}}\\&\le
C'\varepsilon^p\|\theta'\|_{p,w}^2
\(\int_{|\xi|<\varepsilon}N_{\theta,p,w}(\xi)dA(\xi)\)^{\frac{p-2}{p}}
 \,.\end{align*}
 Finally, (\ref{counting-est2}) follows by the same argument.
\end{proof}
There are two remarks to be made here. The first one is the
following interesting corollary.
\begin{cor} Let $(p_1,p_2)\in (1/2,1]\times(1/2,1]
\cup(1,\infty)\times(1,\infty)$, with $0<p_2-p_1<1$, and assume that
the weight $w$ satisfies $b_w>-1,~ a_w<2p_1-2$, when $1/2<p_1<p_2\le
1$, and $p_2-2<b_w\le a_w<p_1-1$ when $p_1,p_2>1$. Then there exists
a constant $c>0$ depending only on $p_1,p_2$ and $w$ such that for
every inner function $\theta$ we have
$$c^{-1}\|\theta'\|_{p_2,w}^{p_2}
\le \int_\D|\theta'(z)|^{p_1}w(z)(1-|z|)^{p_1-p_2}dA(z)\le
c\|\theta'|_{p_2,w}^{p_2}\,.$$
\end{cor}
\begin{proof} If $w_1(z)=w(z)(1-|z|)^{p_1-p_2}$ then $p_1, w_1$
satisfy the conditions in Theorem \ref{inner-est} and
$$N_{\theta,p_1,w_1}(\z)=N_{\theta,p_2,w}(\z)\,.$$
\end{proof}
Note that the inequality on the left follows by the Schwarz-Pick
lemma, but the one on the right is somewhat surprising.
\par
Secondly, we should say a few words about the conditions on $w$ and $p$. The point we want to make is that they place us in an interesting type of spaces. For example, if $w$ is a standard weight, $w(z)= (\a+1) (1-|z|)^\a$, $\a>-1$, then $a_w=b_w=\a$, and the assumptions
are
\[
 \a<2p-2\,,\quad\text{if}~~1/2<p\le 1\,,\quad
 p-2<\a<p-1\,,\quad\text{if}~~p>1\,.
\]
If $\a>p-1,~p>1/2$ then the derivative of an arbitrary $H^2$-function belongs to $L_a^{p,\a}=L_a^p(w)$, and if $2\ge p>1,~\a<p-2$, only the
derivative of a finite Blaschke product can belong to this space.
\section{The use of duality}
 \label{sec-duality}
\par
If $p>1$, the conditions  $p-2<b_w \le a_w<p-1$, considered in the
previous section, ensure that the dual of $L_a^p(w)$ can be identified with $L_a^q(w^{-q/p}),~\frac{1}{p}+\frac1{q}=1$ (see \cite{L}), with respect to the unweighted Bergman pairing
$$
 (f,g)=\int_\D f\overline{g}dA\,,
$$
\textit{i.e.}, the $L_a^p(w)$-norm can be estimated as
\begin{equation}\label{duality-est}
\|f\|_{p,w}\asymp \sup_{\|g\|_{q,w^{-q/p}}\le
1}\left|\int_\D\overline{f}gdA\right|\,.
\end{equation}
  This information proves
useful, for it provides an additional tool for estimating the
$L_a^p(w)$- norm of derivatives of Blaschke products. By a direct
computation based on the Cauchy-Green (Stokes) formula \cite[p.~17]{DS}, we have
\begin{equation}\label{stokes}
 \int_\D\overline{f'}gdA
 =\frac1{2i}\int_{|z|=1}\overline{f}(z)g(z)dz\,,\quad f,g\in H^2\,.
\end{equation}
An immediate application of \eqref{stokes} is the following result
related to the so-called $\mathcal{F}$-property.
\par
\begin{prop}\label{f-prop} If $p>1$, and $p-2<b_w \le
a_w<p-1$, there exists a positive constant $c_{w,p}$ depending only
on $p$ and $w$ such that such that whenever $\theta_1,\theta_2$ are
inner functions with $\theta_2/\theta_1\in H^\infty$, we have
$$\|\theta_1\|_{p,w}\le c_{p,w}\|\theta_2\|_{p,w}\,.$$
\end{prop}
\begin{proof}
By \eqref{stokes} we have that
$$
 \int_\D\overline{\theta_1'}gdA=\int_\D\overline{\theta_2'}
\frac{\theta_1}{\theta_2}gdA\,,
$$
whenever $g\in H^2$. Since $H^2$ is dense in $L_a^p(w)$, the result follows by \eqref{duality-est}.
\end{proof}
If $f=B$ is a finite Blaschke product with simple zeros, then using
the fact that $\overline{B}(z)=1/B(z)$ when $|z|=1$, we can evaluate
the line integral in \eqref{stokes} with help of the residue
formula. Since $H^2$ is dense in $L_a^q(w^{-q/p})$ we have in this
case
\begin{equation}\label{blaschke-stokes}
 \frac1{2i}\int_{|z|=1}\overline{B}(z)g(z)dz=\pi
 \sum_{n}\frac{1-|z_n|^2}{B_n(z_n)}g(z_n)\,,\quad g\in
 L_a^q(w^{-q/p})\,,
\end{equation}
where $(z_n)_{n=1}^\infty$ are the zeros of $B$ and, as before,
$B_n$ denotes the Blaschke product with zeros $\{z_k:~k\neq n\}$. We
will also use the notation $B^{[N]}$ for the $N$-th partial product
of the Blaschke product $B$.
\begin{prop}\label{duality-prop} 
Assume that $p>1$, and $p-2<b_w \le a_w<p-1$. Then: 
\par\noindent 
(i) For every Blaschke product $B$ with simple zeros we
have
$$
 \|B^\prime\|_{p,w}\asymp \limsup_{N\to\infty}\sup_{\|g\|_{q,w^{-q/p}}\le
1}\left|\sum_{n=1}^N\frac{1-|z_n|^2}{B^{[N]}_n(z_n)}g(z_n)\right|\,,
$$
where the constants involved depend only on $p$ and $w$. 
\par\noindent 
(ii) If the zero sequence $(z_n)_{n=1}^\infty$ is interpolating for
$L_a^q(w^{-q/p})$ then
$$
 \|B^\prime\|_{p,w}^p\asymp
 \sum_{n}\frac{(1-|z_n|)^{2-p}w(|z_n|)}{|B_n(z_n)|^p}\,.
$$ 
In particular, this holds whenever $B$ is an interpolating Blaschke
product.
\end{prop}
\begin{proof}
(i) By Proposition \ref{f-prop} it follows immediately that
$$ \limsup_{N\to\infty}\sup_{\|g\|_{q,w^{-q/p}}\le
1}\left|\sum_{n=1}^N\frac{1-|z_n|^2}{B^{[N]}_n(z_n)}g(z_n)\right|\le
C \|B^\prime\|_{p,w}\,.$$ Conversely, let $g$ be a polynomial such that
$\|g\|_{q,w^{-q/p}}\le 1$ and
$$\|B^\prime\|_{p,w}<C'\left|\int_\D\overline{B^\prime}gdA\right|\,,$$
where $C'$ depends only on $p$ and $w$. Since
$$
 \int_\D\overline{B^\prime}gdA=\lim_{N\to\infty}
\int_\D\overline{(B^{[N]})'}gdA\,,
$$
the result follows.
\par
(ii) If
$(z_n)_{n=1}^\infty$ is interpolating for $L_a^q(w^{-q/p})$ then
\begin{align*}\sup_{\|g\|_{q,w^{-q/p}}\le
1}\left|\sum_{n=1}^N\frac{1-|z_n|^2}{B^{[N]}_n(z_n)}g(z_n)\right|&\asymp
\sup_{\|(x_n)\|_{\ell^q}\le
1}\left|\sum_{n=1}^N\frac{(1-|z_n|^2)^{1-2/q}
w^{1/p}(|z_n|)}{B^{[N]}_n(z_n)}x_n\right|\\&=
\sum_{n=1}^N\frac{(1-|z_n|)^{2-p}w(|z_n|)}{|B^{[N]}_n(z_n)|^p}\,,\end{align*}
and the result follows by the monotone convergence theorem.
\end{proof}
It is interesting to note that part (ii) together with Theorem
\ref{main1} (i) yield an average lower bound for $|B_n(z_n)|$ in the
case when $(z_n)_{n=1}^\infty$ is interpolating for $L_a^q(w^{-q/p})$.
We have that
 \begin{equation} \label{cor-interp}
\sum_{n}\frac{(1-|z_n|)^{2-p}w(|z_n|)}{|B_n(z_n)|^p}\le
C\sum_{n}(1-|z_n|)^{2-p}w(|z_n|)|B_n(z_n)|^p\,.
\end{equation}
We now turn to Cohn's theorem \cite{Ch} which asserts that for
interpolating Blaschke products $B$ and $s\in (1/2,2/3)$ we have
$$
 K_B^\prime\subset H^s \Longleftrightarrow B^\prime\in H^{\frac{2s}{2-s}}\,,
$$
where $K_B=H^2\ominus BH^2$, and $K_B^\prime$ denotes the set of
derivatives of functions in $K_B$. Given  a Blaschke product $B$ and
$p>1$ we denote by
$$ 
 K_B^p=H^p\cap B\overline{H_0^p}=(BH^q)^\perp\,,\quad
 \frac{1}{p}+\frac1{q}=1\,,
$$ 
and by $(K_B^p)'$ the set of derivatives of functions in $K_B^p$.
\begin{thm} \label{cohn-an}
Let $B$ be an interpolating Blaschke product, let $p>1$, and let $w$
be a normal weight. Let $1<s<p$ and assume that
$s-2+\frac{s}{p}<b_w\le a_w<s-1$. Then $(K_B^p)'\subset L_a^s(w)$ if
and only if $B^\prime\in L_a^{s\frac{p-1}{p-s}}(w^{\frac{p}{p-s}})$.
\end{thm}
\begin{proof} Since $(z_n)_{n=1}^\infty$ is an interpolating sequence, a function $f\in H^p$ belongs to
$K_B^p$ if and only if it can be written as
$$
 f(z)=\sum_{n}c_n\frac{(1-|z_n|)^{1/q}}{1-\overline{z}_nz}\,,
$$
where the numbers $c_n\in \C$ are uniquely determined by $f$ and
$z_n$, $n=1,2,\ldots$, are the zeros of $B$. Moreover,
$$
 \|f\|_p^p\asymp \sum_n|c_n|^p\,.
$$
Now let $r>1$ be given by $\frac1{r}=\frac1{s}-\frac{1}{p}$ and
denote by $r'$ and $s'$ the conjugate indices of $r$ and $p$, respectively; that is,
$$
 \frac{1}{r'}+\frac1{r}=\frac1{s'}+\frac1{s}=1\,.
$$ 
By the condition on $w$, we see that we can apply \eqref{duality-est} to obtain
\begin{align*}\sup\{&\|f'\|_{s,w}\,\colon\, f\in K_b^p,~\|f\|_{H^p}\le 1\}\\& \asymp
\sup\left\{\left|\int_\D f'
\overline{g}dA\right|\,\colon\,\|g\|_{s',w^{-s'/s}}\le1, f\in K_b^p,~
\|f\|_{H^p}\le 1 \right\}\\&\asymp \sup\left\{\left|\sum_n
\overline{z}_nc_n\overline{g(z_n)} (1-|z_n|)^{1/q}
\right|\,\colon\,\|g\|_{s',w^{-s'/s}}\le1, \sum_n|c_n|^p\le 1
\right\}\,.\end{align*} Now recall that, by Theorem
\ref{interpolation}, $(z_n)_{n=1}^\infty$ is interpolating for
$L_a^{s'}(w^{-s'/s})$, hence
\begin{align*}
 &\sup\left\{\left|\sum_n \overline{z}_nc_n\overline{g(z_n)}
 (1-|z_n|)^{1/q} \right|\,\colon\,\|g\|_{s',w^{-s'/s}}\le 1\,, \   \sum_n|c_n|^p\le 1 \right\} \\
 &\asymp \sup\left\{\left|\sum_na_nc_nw(|z_n|)^{1/s}
 (1-|z_n|)^{1/q-2/s'} \right|\,\colon\,\sum_n|a_n|^{s'}\le 1\,, \  \sum_n|c_n|^p\le 1 \right\}\,.
\end{align*}
With our notations, we can rewrite this as
$$\sup\{\|f'\|^r_{s,w}\,\colon\, f\in K_b^p,~\|f\|_{H^p}\le 1\}\asymp
\sum_nw(|z_n|)^{r/s} (1-|z_n|)^{r/q-2r/s'}\,,$$ and since
$\frac1{q}-\frac1{s'}=\frac{1}{r}$, we obtain
$$\sup\{\|f'\|^r_{s,w}\,\colon\, f\in K_b^p,~\|f\|_{H^p}\le 1\}\asymp
\sum_nw(|z_n|)^{r/s} (1-|z_n|)^{2-r/q}\,.$$ Finally, from the fact
that
$$a_{w^{r/s}}=\frac{r}{s}a_w\,,\quad b_{w^{r/s}}=\frac{r}{s}b_w\,,$$
 we see that the conditions in the statement are equivalent to
 $$\frac{r}{q}-2<b_{w^{r/s}}\le a_{w^{r(s}}<\frac{r}{q}-1\,,$$
 and since $\frac{r}{q}>1$,  Theorem \ref{main1} yields
 $$\sum_nw(|z_n|)^{r/s} (1-|z_n|)^{2-r/q}\asymp
 \int_\D|B^\prime|^{r/q}w^{r/s}dA\,,$$
 which completes the proof.
\end{proof}
We should point out here that the case when $s\ge p$ can be treated
similarly, but it leads either to the limit cases in Theorem
\ref{main1}, or to trivial results.
\section{Derivatives of higher order}
 \label{sec-higher}
\par
It is well known that
\[
 \int_\D |f(z)|^p dA(z) \asymp \sum_{k=0}^{n-1} |f^{(k)}(0)|^p +
 \int_\D |f^{(n)}(z)|^p (1-|z|^2)^{n p} dA(z)
\]
and it turns out that this extends to weighted Bergman spaces with
normal functions as well.
\par
\begin{other} \label{other-compar}
If $b_w<1$ then
\[
 \|f\|^p_{p,w}\asymp \sum_{k=0}^{n-1} |f^{(k)}(0)|^p +
 \int_\D |f^{(n)}(z)|^p (1-|z|^2)^{n p} w(|z|)\,d A(z) \,.
\]
\end{other}
In the unweighted case $w=1$, the above statement is a well known
theorem of Hardy and Littlewood (\textit{cf.\/} \cite[Chapter~5]{D}
where one inequality is proved, and the other can easily be proved
using the techniques displayed there). The above generalization was
proved explicitly in \cite{AS} for $p\ge 1$ and in \cite{H} for
$p<1$.
\par
Taking into account that the radial function $w$ is normal if and only if the associated function $W(z)=(1-|z|)^p w(z)$ is normal, one can derive from Theorem~\ref{other-compar} various analogues of our main theorems for derivatives of higher order. Just to illustrate a flavor of such results, we now state without proof one such statement.
\begin{thm}\label{main-higher}
Let $w$ be a normal function and let $1/2<p<\infty$.
\begin{enumerate}
 \item[(i)] Assume that   $a_w<2 p-2,~b_w>-1$, if  $1/2<p\le 1$,
 and that $a_w<p-1,~b_w>p-2$, if  $p>1$.  Then there
 exists a positive constant $c_{p,w}$ such that for every Blaschke
 product $B$ we have
\[
  \sum_{k=0}^{n-1} |B^{(k)}(0)|^p + \|B^{(n)}\|^p_{p,w}\le c_{p,w} \sum_{B(z)=0} (1-|z|)^{2-np} \,w(|z|)\,.
\]
\item[(ii)] Suppose that the zero set  of the
Blaschke product $B$ is separated with separation constant
$\delta>0$. Suppose that either $a_w<2p-2$, $b_w>-1$, when $1/2<p\le 1$, or $a_w<p-2$, $b_w>-1$, when $p>1$. Then, in both cases, there exists a positive constant $c_{p,w,\d}$ such that
\[
  \sum_{B(z)=0}^\infty (1-|z|)^{2-np}
 \,w(|z|)\le \sum_{k=0}^{n-1} |B^{(k)}(0)|^p + c_{p,w,\delta}\|B^{(n)}\|^p_{p,w}\,.
\]
 \end{enumerate}
\end{thm}


\begin{thebibliography}{99}

\bibitem{A1}
P.~Ahern, The mean modulus of the derivative of an inner function,
\emph{Indiana Univ. Math. J.} \textbf{28} (1979), no. 2, 311--347.

\bibitem{A2}
P.~Ahern, The Poisson integral of a singular measure, \emph{Canad.
J. Math.} \textbf{35} (1983), no. 4, 735--749.

\bibitem{AC1}
P. Ahern and D. Clark, On inner functions with $H\sp{p} $
derivative, \emph{Michigan Math. J.} \textbf{21} (1974), 115--127.

\bibitem{AC2}
P. Ahern and D. Clark, On inner functions with $B\sp{p} $
derivative, \emph{Michigan Math. J.} \textbf{23} (1976), no. 2,
107--118.

\bibitem{AS} A. Aleman and A. Siskakis, Integration operators on
Bergman spaces, \textit{Indiana Univ. Math. J.\/} \textbf{46}
(1997), no. 2, 337--356.

\bibitem{Be} D. B\'ekoll\'e, In\'egalit\'e \`{a} poids pour le
projecteur de Bergman dans la boule unit\'e de $ C\sp{n}$,
\emph{Studia Math.}, {\bf 71} (1981/82), no. 3, 305--323.

\bibitem{Ch}
W.S.~Cohn, On the $H^p $ classes of derivative of functions
orthogonal to invariant subspaces, \emph{Michigan Math. J.\/}
\textbf{30} (1983), 221--229.

\bibitem{Cl}
P. Colwell, \emph{Blaschke products. Bounded analytic functions\/},
University of Michigan Press, Ann Arbor, Michigan 1985.

\bibitem{D}
P.L. Duren, \emph{Theory of $H\sp{p} $ Spaces\/}, Academic Press,
New~York-London 1970. Reprint: Dover, Mineola, New York 2000.

\bibitem{DS}
P.L.~Duren and A.P.~Schuster, \emph{Bergman Spaces\/}, Mathematical
surveys and monographs; no.100, American Mathematical Society,
Providence, RI (2004).

\bibitem{FM} E. Fricain, J. Mashreghi, Integral means of the
derivatives of Blaschke products.  \textit{Glasg. Math. J.} \textbf{50} (2008), no. 2, 233--249.

\bibitem{G}
J.B. Garnett, \emph{Bounded Analytic Functions}, Academic Press, New
York, etc. 1981.

\bibitem{GP1}
D. Girela and J.A. Pel\'aez, On the derivative of infinite Blaschke
products, \emph{Illinois J. Math.} \textbf{48} (2004), No.1,
121-130.

\bibitem{GP2}
D. Girela and J. A. Pel\'aez, On the membership in Bergman spaces of
the derivative of a Blaschke product with zeros in a Stolz domain,
\textit{Canadian Math. Bull.} \textbf{49} (2006),  no. 3, 381--388.

\bibitem{GPV1}
D. Girela, J.A. Pel\'aez, and D. Vukoti\'c, Integrability of the
derivative of a Blaschke product, \textit{Proc. Edinburgh Math. Soc.
(2)\/}, \textbf{50} (2007), no. 3, 673-687.

\bibitem{GPV3}
D. Girela, J.A. Pel\'aez, and D. Vukoti\'c, Uniformly discrete
sequences in regions with tangential approach to the unit circle,
\textit{Complex Variables Elliptic Ecu.\/} \textbf{52}  (2007),
No's. 2--3, 161--173.

\bibitem{GPV2}
D. Girela, J.A. Pel\'aez, and D. Vukoti\'c, Interpolating Blaschke
products: Stolz and tangential approach regions, \textit{Constr.
Approx.\/}, \textbf{27} (2008), no. 2, 203--216.

\bibitem{HKZ}
H.~Hedenmalm, B.~Korenblum, and K.~Zhu, \emph{Theory of Bergman
Spaces\/}, Graduate Texts in Mathematics, Vol. \textbf{199},
Springer, New York, Berlin, etc. 2000.

\bibitem{H}
Z. Hu, Extended Ces\`aro operators on mixed norm spaces,
\textit{Proc. Amer. Math. Soc.\/} \textbf{131} (2003), no. 7,
2171--2179.

\bibitem{Ki}
H.O. Kim, Derivatives of Blaschke products, \emph{Pacific J. Math.}
\textbf{114} (1984), 175--190.

\bibitem{Ku}
M.A. Kutbi, Integral means for the first derivative of Blaschke
products, \emph{Kodai Math. J.} \textbf{24} (2001), no. 1, 86--97.

\bibitem{L}
D. H. Luecking, Representation and duality in weighted spaces of
analytic functions, \textit{Indiana Univ. Math. J.} \textbf{34}
(1985), no. 2, 319--336.

\bibitem{N}
C. Nolder, An $L^p$ definition of interpolation Blaschke products, \textit{Proc. Amer. Math. Soc.\/} \textbf{128}, no. 6 (2000), 1799--1806.

\bibitem{Pe}
J. A. Pel\'aez, Sharp results on the integrability of the derivative of an interpolating Blaschke product, \textit{Forum Math.} \textbf{20}   (2008), no. 6, 1039--1054.

\bibitem{Pi}
G.~Piranian, Bounded functions with large circular variation,
\emph{Proc. Amer. Math. Soc.} \textbf{19}, no. 6 (1968), 1255--1257.

\bibitem{Pr}
D. Protas, Blaschke products with derivative in $H\sp{p} $ and
$B\sp{p} $, \emph{Michigan Math. J.} \textbf{20} (1973), 393--396.

\bibitem{R}
W. Rudin, The radial variation of analytic functions, \emph{Duke
Math. J.} \textbf{22} (1955), 235--242.

\bibitem{S}
K. Seip, Beurling type density theorems in the unit disk,
\emph{Invent. Math\/}. \textbf{113} (1993), 21--39.

\bibitem{Seip} K. Seip, \emph{Interpolation and sampling in spaces
 of analytic functions}, University
Lecture Series, \textbf{33}. American Mathematical Society,
Providence, RI, 2004.

\bibitem{Sh}
J.H.~Shapiro, \emph{Composition Operators and Classical
Function Theory\/}, Springer-Verlag, New York 1993.

\bibitem{SW}
A. L. Shields and D. L. Williams, Bounded projections, duality, and
multipliers in spaces of analytic functions, \textit{Trans. Amer.
Math. Soc.\/} \textbf{162} (1971), 287--302.

\bibitem{Z}
K.~Zhu, \emph{Operator Theory in Function Spaces\/}, Marcel Dekker,
Inc., Pure and Applied Mathematics \textbf{139}, New York and Basel
1990.


\end{thebibliography}
\end{document}